\documentclass{amsart}
\usepackage{amsfonts}

\theoremstyle{plain}
\newtheorem{theorem}{Theorem}
\newtheorem{lemma}{Lemma}
\newtheorem{corollary}{Corollary}
\newtheorem{proposition}{Proposition}

\theoremstyle{definition}

\theoremstyle{remark}
\newtheorem{remark}{Remark}

\numberwithin{equation}{section}
\input{tcilatex}

\begin{document}
\title[Poisson-Mehler kernels]{On summable, positive Poisson-Mehler kernels
built of Al-Salam--Chihara and related polynomials}
\author{Pawe\l\ J. Szab\l owski}
\address{Department of Mathematics and Information Sciences,\\
Warsaw University of Technology\\
pl. Politechniki 1, 00-661 Warsaw, Poland }
\email{pawel.szablowski@gmail.com}
\date{February 10, 2011}
\thanks{Author is very grateful to an unknown referee for his detailed
remarks and corrections including pointing out numerous misprints that all
together helped to improve the paper.}

\begin{abstract}
Using special technique of expanding ratio of densities in an infinite
series of polynomials orthogonal with respect to one of the densities, we
obtain simple, closed forms of certain kernels built of the so called
Al-Salam---Chihara (ASC) polynomials. We consider also kernels built of some
other families of polynomials such as the so called big continuous $q-$%
Hermite polynomials that are related to the ASC polynomials. The constructed
kernels are symmetric and asymmetric. Being the ratios of the densities they
are automatically positive. We expand also reciprocals of some of the
kernels, getting nice identities built of the ASC polynomials involving 6
variables like e.g. formula (\ref{nice}). These expansions lead to
asymmetric, positive and summable kernels. The particular cases (referring
to $q=1$ and $q=0$) lead to the kernels build of certain linear combinations
of the ordinary Hermite and Chebyshev polynomials.
\end{abstract}

\maketitle

\section{Introduction}

In many models of the so called $q-$oscillators considered in quantum
physics, classical and noncommutative probability or generally in some
branches of analysis appears a problem of summing and examining positivity
of kernels built of certain families of orthogonal polynomials (see e.g. 
\cite{Floreanini97}, \cite{Flo97}, \cite{Askey1993}). The kernels (more
precisely the Poisson--Mehler kernels) are, generally speaking, expressions
of the form $K\left( x,y\right) \allowbreak =\allowbreak \sum_{n\geq
0}a_{n}D_{n}\left( x\right) F_{n}\left( y\right) ,$ where $\left\{
D_{n}\right\} _{n\geq 0},$ $\left\{ F_{n}\right\} _{n\geq 0}$ are certain
families of orthogonal polynomials and $x,\allowbreak y,\allowbreak \left\{
a_{n}\right\} _{n\geq 0}$ are real numbers. Usually the numbers $a_{n}$ are
of the form $t^{n}/\left\Vert D_{n}\right\Vert \left\Vert F_{n}\right\Vert ,$
$\left\vert t\right\vert <1$ where $\left\Vert .\right\Vert $ denotes
certain (usually$L_{2}$) norm. Often $D_{n}\allowbreak =\allowbreak F_{n},$
then such kernels are called symmetric otherwise they are called asymmetric.
Sometimes the symmetric kernels, with $a_{n}\allowbreak =\allowbreak
t^{n}/\left\Vert D_{n}\right\Vert ^{2}$ are called bilinear generating
functions.

Notice that we have $\int K\left( x,y\right) D_{n}\left( x\right) d\delta
\left( x\right) \allowbreak =\allowbreak a_{n}\left\Vert D_{n}\right\Vert
^{2}F_{n}\left( y\right) $ where $\delta \left( x\right) $ denotes the
measure with respect to which polynomials $\left\{ D_{n}\right\} $ are
orthogonal. The property often exploited in various branches of analysis
including integral equations.

One would like to express such kernels in closed, compact forms and what is
more important give conditions under which such kernels are nonnegative for
all $x,y\allowbreak $ from certain Cartesian product of intervals.

Then for example an expression $h(x,y)\allowbreak =\allowbreak f\left(
x\right) K\left( x,y\right) g\left( y\right) $ with $K$ being a positive
kernel, $f$ and $g$ being probability densities such that the polynomials $%
\left\{ D_{n}\right\} $ and $\left\{ F_{n}\right\} $ are orthogonal with
respect to the measures with these densities would give a Lancaster
expansion of certain $2$ dimensional probabilistic density. This so since $%
\int \int h\left( x,y\right) dxdy\allowbreak =\allowbreak 1$ and $h\left(
x,y\right) \geq 0.$ For the Lancaster expansions and their most important
properties see e.g. \cite{Lancaster}.

The problem of summing and examining positivity of kernels is nontrivial and
was solved only for a few families of polynomials including the $q$-Hermite
polynomials. Thus it is important to extend the list of polynomials for
which kernels are summable and nonnegative.

This paper extends such a list by adding to it two families of orthogonal
polynomials (indexed by real parameter $q$). One family added partially are
the so called big $q-$Hermite (bH) polynomials and the other family more
completely added are the so called Al-Salam--Chihara (ASC) polynomials and
of course the $q\allowbreak =\allowbreak 1$ and $\allowbreak q\allowbreak
=\allowbreak 0$ special cases. These special cases simplify for $%
q\allowbreak =\allowbreak 1$ to respectively shifted and shifted and
re-scaled Hermite polynomials and for $q\allowbreak =\allowbreak 0$ to
certain combinations of the Chebyshev polynomials of the second kind..

The summation formulae for kernels built of bH and ASC polynomials have been
found in 1996 in \cite{suslov96}. They were expressed with the help of the
basic hypergeometric functions: $_{3}\phi _{2}$ in the case of bH
polynomials and $_{8}\phi _{7}$ (in fact $_{8}W_{7})$ in the case of ASC
polynomials under some conditions imposed on some parameters of considered
families of polynomials. Hence they are 'inconsumable' for the reader not
familiar with specialized branches of special functions theory. The
positivity of obtained kernels was shown in \cite{suslov96} only for the
very special cases.

We will present a compact, simple, summable form of those kernels for
different (than in \cite{suslov96}) sets of restrictions imposed on the
parameters and for special values of the parameter $t.$ In summing kernels
we apply the original method of expanding the ratio of densities in a
certain Fourier series of polynomials that are orthogonal with respect one
of the densities. This technique was used successfully in \cite{Expansion}
to present new proofs of compact forms of some kernels. Thus we obtain some
new results and also known results by the new technique. Besides we obtain
expansions of reciprocals of some of the kernels. This procedure leads to
non-symmetric kernels that are naturally nonnegative.

To help those who are not familiar with the $q-$series theory we will
present basic notation and basic notions used in this theory and also point
out and discuss in detail two special cases. As stated above the case $%
q\allowbreak =\allowbreak 1$ leads to shifted and re-scaled Hermite
polynomials and $q\allowbreak =\allowbreak 0$ leads to linear combinations
of Chebyshev polynomials of the second kind. Hence, although focused mostly
on the applications of $q-$series theory the paper can also be interesting
for the people working in orthogonal polynomials.

The paper is organized as follows. In the next section we present basic
notation and notions of the $q-$series theory, then we introduce and recall
basic properties of families of polynomials discussed in the paper. Next
section presents main results while the last section contains some longer,
less interesting proofs.

\section{Preliminaries}

In order to formulate briefly the properties of bH and ASC polynomials we
must introduce notation traditionally used in the so called $q-$series
theory. $q$ is a real parameter $-1<q\leq 1.$ Having $q$ we define $\left[ 0%
\right] _{q}\allowbreak =\allowbreak 0;$ $\left[ n\right] _{q}\allowbreak
=\allowbreak 1+q+\ldots +q^{n-1}\allowbreak ,$ $\left[ n\right]
_{q}!\allowbreak =\allowbreak \prod_{i=1}^{n}\left[ i\right] _{q},$ with $%
\left[ 0\right] _{q}!\allowbreak =1,\QATOPD[ ] {n}{k}_{q}\allowbreak
=\allowbreak \left\{ 
\begin{array}{ccc}
\frac{\left[ n\right] _{q}!}{\left[ n-k\right] _{q}!\left[ k\right] _{q}!} & 
, & n\geq k\geq 0 \\ 
0 & , & otherwise%
\end{array}%
\right. $. It will be useful to use the so called $q-$Pochhammer symbol for $%
n\geq 1:$%
\begin{equation*}
\left( a;q\right) _{n}=\prod_{i=0}^{n-1}\left( 1-aq^{i}\right) ,\left(
a_{1},a_{2},\ldots ,a_{k};q\right) _{n}\allowbreak =\allowbreak
\prod_{i=1}^{k}\left( a_{i};q\right) _{n},
\end{equation*}%
with $\left( a;q\right) _{0}=1.$

Often $\left( a;q\right) _{n}$ as well as $\left( a_{1},a_{2},\ldots
,a_{k};q\right) _{n}$ will be abbreviated to $\left( a\right) _{n}$ and $%
\left( a_{1},a_{2},\ldots ,a_{k}\right) _{n},$ if it will not cause
misunderstanding.

It is easy to notice that $\left( q\right) _{n}=\left( 1-q\right) ^{n}\left[
n\right] _{q}!$ and thus that:%
\begin{equation*}
\QATOPD[ ] {n}{k}_{q}\allowbreak =\allowbreak \allowbreak \left\{ 
\begin{array}{ccc}
\frac{\left( q\right) _{n}}{\left( q\right) _{n-k}\left( q\right) _{k}} & ,
& n\geq k\geq 0 \\ 
0 & , & otherwise%
\end{array}%
\right. .
\end{equation*}%
\newline
Notice that for $n\allowbreak \geq \allowbreak k\allowbreak \geq \allowbreak
0$ we have: $\left[ n\right] _{1}\allowbreak =\allowbreak n,\left[ n\right]
_{1}!\allowbreak =\allowbreak n!,$ $\QATOPD[ ] {n}{k}_{1}\allowbreak
=\allowbreak \binom{n}{k},$ $\left( a;1\right) _{n}\allowbreak =\allowbreak
\left( 1-a\right) ^{n}$ and $\left[ n\right] _{0}\allowbreak =\allowbreak
\left\{ 
\begin{array}{ccc}
1 & if & n\geq 1 \\ 
0 & if & n=0%
\end{array}%
\right. ,$ $\left[ n\right] _{0}!\allowbreak =\allowbreak 1,$ $\QATOPD[ ] {n%
}{k}_{0}\allowbreak =\allowbreak 1,$ $\left( a;0\right) _{n}\allowbreak
=\allowbreak \left\{ 
\begin{array}{ccc}
1 & if & n=0 \\ 
1-a & if & n\geq 1%
\end{array}%
\right. .$

The continuous $q-$Hermite polynomials $h_{n}\left( x|q\right) $ are defined
by the following $3-$term recurrence 
\begin{equation*}
h_{n+1}(x|q)=2xh_{n}\left( x|q\right) -(1-q^{n})h_{n-1}\left( x|q\right) ,
\end{equation*}%
with $h_{-1}\left( x|q\right) \allowbreak =\allowbreak 0,$ $h_{0}\left(
x|q\right) \allowbreak =\allowbreak 1,$ for $\left\vert q\right\vert $%
\allowbreak $<1.$ Many facts are known about these polynomials. The most
important recent references are \cite{Koek} and \cite{IA}.

We will consider the modified version of these polynomials mostly because
these modifications have nice probabilistic interpretations and applications
see e.g. \cite{Bryc2001S}, \cite{Bo}, \cite{bms}, \cite{AW}, \cite{Szab5}, 
\cite{Expansion}. Namely we will consider the polynomials $\left\{
H_{n}\left( x|q\right) \right\} $ defined by the relationship for $q\neq 1$
: 
\begin{equation*}
H_{n}\left( x|q\right) \allowbreak =\allowbreak h_{n}\left( x\sqrt{1-q}%
/2|q\right) /\left( 1-q\right) ^{n/2};n\geq 0
\end{equation*}%
and $H_{n}\left( x|1\right) \allowbreak =\allowbreak H_{n}\left( x\right) ,$
where $\left\{ H_{n}\left( x\right) \right\} $ are the so called
probabilistic Hermite polynomials that is polynomials orthogonal with
respect to the measure that has density equal to $\exp \left(
-x^{2}/2\right) /\sqrt{2\pi }$.

To see that really $H_{n}\left( x|q\right) \allowbreak =\allowbreak $ $%
H_{n}\left( x\right) $ one has to refer to $3-$term recurrence satisfied by
the polynomials $\left\{ H_{n}\left( x|q\right) \right\} $ that in this case
has the following simple form for $n\geq 0:$%
\begin{equation}
H_{n+1}\left( x|q\right) =xH_{n}\left( x|q\right) -\left[ n\right]
_{q}H_{n-1}\left( x|q\right) ,  \label{_H}
\end{equation}%
with $H_{-1}\left( x|q\right) =0,$ $H_{0}\left( x|q\right) \allowbreak
=\allowbreak 1.$ Symbol $\left[ n\right] _{q}$ was explained above and also
above we remarked that $\left[ n\right] _{1}\allowbreak =\allowbreak n$
hence for $q\allowbreak =\allowbreak 1$ (\ref{_H}) reduces to well known
(see \cite{IA} or \cite{Andrews1999}) $3-$term recurrence satisfied by
'probabilistic' Hermite polynomials i.e. the ones that are orthogonal with
respect to measure with the density $\exp \left( -x^{2}/2\right) /\sqrt{2\pi 
}$.

From (\ref{_H}) it follows also that $H_{n}\left( x|0\right) \allowbreak
=\allowbreak U_{n}\left( x/2\right) ,$ where $\left\{ U_{n}\left( x\right)
\right\} _{n\geq -1}$ are the so called Chebyshev polynomials of the second
kind. It is so since $\left[ n\right] _{0}\allowbreak =\allowbreak 1$ for $%
n\geq 1.$ More precisely polynomials $\left\{ U_{n}\right\} $ satisfy the
following 3-term recurrence: 
\begin{equation*}
U_{n+1}\left( x\right) \allowbreak =\allowbreak 2xU_{n}\left( x\right)
-U_{n-1}\left( x\right) ,
\end{equation*}%
with $U_{-1}\left( x\right) =0,$ $U_{0}\left( x\right) =1.$ To learn more
about Chebyshev polynomials the reader is referred either to \cite%
{Andrews1999} or to \cite{IA}.

It is known in particular that the characteristic function of the
polynomials $\left\{ H_{n}\right\} $ for $\left\vert q\right\vert <1$ is
given by the formula:%
\begin{equation*}
\sum_{n=0}^{\infty }\frac{t^{n}}{\left[ n\right] _{q}!}H_{n}\left(
x|q\right) =\varphi \left( x|t,q\right) ,
\end{equation*}%
where 
\begin{equation}
\varphi \left( x|t,q\right) \allowbreak =\allowbreak \frac{1}{%
\prod_{k=0}^{\infty }\left( 1-(1-q)xtq^{k}+(1-q)t^{2}q^{2k}\right) }.
\label{chH}
\end{equation}%
Convergence here is absolute for $\left\vert \sqrt{1-q}t\right\vert <1.$ For
the proof of this formula see e.g. \cite{IA}. (\ref{chH}) will be justified
as a simple corollary of our main result presented below. Besides because of
the fact that $H_{n}\left( x|1\right) \allowbreak =\allowbreak H_{n}\left(
x\right) $ we set 
\begin{equation*}
\varphi \left( x|t,1\right) \allowbreak =\allowbreak \exp \left(
xt-t^{2}/2\right) .
\end{equation*}%
Notice also that: 
\begin{equation*}
\varphi \left( x|t,0\right) \allowbreak =\allowbreak \frac{1}{1-xt+t^{2}}.
\end{equation*}

It is also known (again see \cite{IA}) that for $-1<q<1,$ these polynomials
are orthogonal with respect to the positive measure defined by the following
density:%
\begin{equation}
f_{N}\left( x|q\right) =\frac{\sqrt{1-q}\left( q\right) _{\infty }}{2\pi 
\sqrt{4-(1-q)x^{2}}}\prod_{k=0}^{\infty }\left(
(1+q^{k})^{2}-(1-q)x^{2}q^{k}\right) I_{S\left( q\right) }\left( x\right) ,
\label{fN}
\end{equation}%
$x\in \mathbb{R},$ where $S\left( q\right) =[-2/\sqrt{1-q},2/\sqrt{1-q}]$
and $I_{A}\left( x\right) \allowbreak =\allowbreak \left\{ 
\begin{array}{ccc}
1 & if & x\in A \\ 
0 & if & x\notin A%
\end{array}%
\right. $. That is 
\begin{equation*}
\int_{S\left( q\right) }H_{n}\left( x|q\right) H_{m}\left( x|q\right)
f_{N}\left( x|q\right) dx\allowbreak =\allowbreak \left[ n\right]
_{q}!\delta _{m,n},
\end{equation*}%
where $\delta _{n,m}\allowbreak =\allowbreak 0$ if $n\neq m$ and $1$ for $%
m=n.$ The case $q\allowbreak =\allowbreak 1$ was presented above. We set
then $f_{N}\left( x|1\right) \allowbreak =\allowbreak \exp \left(
-x^{2}/2\right) /\sqrt{2\pi }.$ One can also easily notice that 
\begin{equation*}
f_{N}\left( x|0\right) \allowbreak =\allowbreak \frac{1}{2\pi }\sqrt{4-x^{2}}%
I_{[-2,2]}\left( x\right) .
\end{equation*}

The big continuous $q-$Hermite polynomials (briefly continuous bH
polynomials) are defined as the polynomials satisfying the following $3$%
-term recurrence (see e.g. \cite{Koek}):%
\begin{equation}
2xh_{n}(x|a,q)=h_{n+1}(x|a,q)+aq^{n}h_{n}(x|a,q)+(1-q^{n})h_{n-1}(x|a,q).
\label{3trh}
\end{equation}%
One can immediately notice that for $a\allowbreak =\allowbreak 0$ we have $%
h_{n}(x|a,q)\allowbreak =\allowbreak h_{n}\left( x|q\right) .$ One proves
(see e.g. \cite{Floreanini97}) that%
\begin{equation}
h_{n}(x|a,q)\allowbreak =\allowbreak \sum_{i=0}^{n}\QATOPD[ ] {n}{i}%
_{q}(-1)^{i}q^{\binom{i}{2}}a^{i}h_{n-i}\left( x|q\right) ,\text{and }%
h_{n}\left( x|q\right) \allowbreak =\allowbreak \sum_{i=0}^{n}\QATOPD[ ] {n}{%
i}_{q}a^{i}h_{n-i}(x|a,q),  \label{bHm}
\end{equation}

where $\left\{ h_{i}\left( x|q\right) \right\} $ denote the continuous $q-$%
Hermite polynomials.

Let us define re-scaled continuous bH polynomials by defining for $n\geq 0:$ 
\begin{equation}
H_{n}(x|a,q)\allowbreak =\allowbreak h_{n}\left( x\sqrt{1-q}/2|a\sqrt{1-q}%
|q\right) /(1-q)^{n/2}.  \label{bH}
\end{equation}%
We will call them the big $q-$Hermite (briefly bH) polynomials.

We immediately have:

\begin{proposition}
\begin{eqnarray}
H_{n}(x|a,q)\allowbreak &=&\sum_{i=0}^{n}\QATOPD[ ] {n}{i}_{q}q^{\binom{i}{2}%
}(-a)^{i}H_{n-i}\left( x|q\right) ,  \label{bigH} \\
H_{n}\left( x|q\right) &=&\sum_{i=0}^{n}\QATOPD[ ] {n}{i}%
_{q}a^{i}H_{n-i}(x|a,q)  \label{bigH2}
\end{eqnarray}
\end{proposition}

\begin{proof}
It is trivial in view of (\ref{bHm}).
\end{proof}

One can also rewrite (\ref{3trh}) in terms of the polynomials $H_{n}(x|a,q):$%
\begin{equation}
xH_{n}(x|a,q)=H_{n+1}(x|a,q)+aq^{n}H_{n}(x|a,q)+[n]_{q}H_{n-1}(x|a,q).
\label{3trH}
\end{equation}

In \cite{Koek} (formula 3.18.2) we can also find the density of the measure
with respect to which polynomials $h_{n}(x|a,q)$ are orthogonal. One can
easily rewrite this density for the polynomials $H_{n}(x|a,q).$ We have then:%
\begin{equation}
f_{bN}(x|a,q)\allowbreak =\allowbreak f_{N}\left( x|q\right) \varphi (x|a,q).
\label{fbN}
\end{equation}%
Notice that here we have to assume $\left\vert a\sqrt{1-q}\right\vert \,<1,$
since then the series $\sum_{n=0}^{\infty }\frac{a^{n}}{\left[ n\right] _{q}!%
}H_{n}\left( x|q\right) $ is absolutely convergent. By the way one shows
(see \cite{Koek}) that for $\left\vert a\sqrt{1-q}\right\vert >1$
polynomials $H_{n}$ are orthogonal with respect to a mixed measure. That is
the measure that apart from absolutely continuous part has also several
(depending on the size of parameter $a)$ discrete "mass" points. However
since we are interested in the expansions of ratios of the densities we will
consider the absolutely continuous case only.

The cases $q\allowbreak =\allowbreak 0$ and $q\allowbreak =\allowbreak 1$
are treated in the following Remark.

\begin{remark}
\label{bH_special}i) $\forall n\geq 0:H_{n}\left( x|a,1\right) \allowbreak
=\allowbreak H_{n}\left( x-a\right) $

ii)$\forall n\geq 0:$ $H_{n}\left( x|a,0\right) \allowbreak =\allowbreak
U_{n}\left( x/2\right) -aU_{n-1}\left( x/2\right) $
\end{remark}

\begin{proof}
i) Follows (\ref{bigH}) and the well known formula concerning Hermite
polynomials: $H_{n}\left( x-a\right) \allowbreak =\allowbreak \sum_{i=0}^{n}%
\binom{n}{i}(-a)^{i}H_{n-i}\left( x\right) .$

ii) Follows (\ref{bigH}) and the fact that $q^{\binom{i}{2}}$ for $%
q\allowbreak =\allowbreak 0$ is nonzero only for $i\allowbreak =\allowbreak
0 $ and $i\allowbreak =\allowbreak 1.$
\end{proof}

Now let us introduce the so called Al-Salam--Chihara (ASC) polynomials $%
P_{n}\left( x|y,\rho ,q\right) $ that are defined by the following 3-term
recurrence:%
\begin{equation}
P_{n+1}\left( x|y,\rho ,q\right) =(x-\rho yq^{n})P_{n}(x|y,\rho
,q)-[n]_{q}(1-\rho ^{2}q^{n-1})P_{n-1}\left( x|y,\rho ,q\right) ,
\label{ASC}
\end{equation}%
with $P_{-1}\left( x|y,\rho ,q\right) =0,$ $P_{0}\left( x|y,\rho ,q\right)
=1.$ In fact this is again modified version of the other family of
polynomials. Namely in the literature (see e.g. \cite{Koek} or \cite{IA})
better known (under the name of ASC) are polynomials defined by the
following 3-term recurrence:%
\begin{equation}
p_{n+1}\left( x|a,b,q\right) =(2x-(a+b)q^{n})p_{n}\left( x|a,b,q\right)
-(1-q^{n})(1-abq^{n-1})p_{n-1}\left( x|a,b,q\right) ,  \label{pn}
\end{equation}%
with $p_{-1}\left( x|a,b,q\right) \allowbreak =\allowbreak 0,$ $p_{0}\left(
x|a,b,q\right) \allowbreak =\allowbreak 1.$ One can easily deduce that if we
take $a$ and $b$ defined by $ab\allowbreak =\allowbreak \rho ^{2}$ and $%
a+b\allowbreak =\allowbreak \rho y\sqrt{1-q},$ then: 
\begin{equation}
P_{n}\left( x|y,\rho ,q\right) =p_{n}\left( x\sqrt{1-q}/2|a,b,q\right)
/(1-q)^{n/2}.  \label{Pn}
\end{equation}%
We will assume that $y\in S\left( q\right) $ and $\left\vert \rho
\right\vert <1.$ Notice that then the parameters $a,b$ are complex but
forming the conjugate pair and also that such a choice of parameters $a,b$
is reasonable since from the Favard's theorem (e.g. \cite{IA}) it follows
that for $\left\vert ab\right\vert \leq 1$ the measure with respect to which
the polynomials $p_{n}$ are orthogonal is positive. Again polynomials $P_{n}$
have nice probabilistic interpretation (again see for example \cite{bms}, 
\cite{AW}, \cite{Szab5}, \cite{Szab-OU-W}, \cite{Szab-bAW}) that is why we
will use them.

Again we have a remark concerning special cases $q=0$ and $q=1.$

\begin{remark}
\label{special}i) $\forall n\geq 0:$ $P_{n}\left( x|y,\rho ,1\right)
\allowbreak =\allowbreak H_{n}\left( \frac{x-\rho y}{\sqrt{1-\rho ^{2}}}%
\right) \left( 1-\rho ^{2}\right) ^{n/2}.$

ii) $\forall n\geq 1:$ $P_{n}\left( x|y,\rho ,0\right) \allowbreak
=\allowbreak U_{n}\left( x/2\right) -\rho yU_{n-1}\left( x/2\right) +\rho
^{2}U_{n-2}\left( x/2\right) ,$ with $P_{0}\left( x|y,\rho ,0\right)
=\allowbreak 1.$
\end{remark}

\begin{proof}
i) For $q=1$ let us denote 
\begin{equation*}
\tilde{P}_{n}(x)=H_{n}\left( \frac{x-\rho y}{\sqrt{1-\rho ^{2}}}\right)
\left( 1-\rho ^{2}\right) ^{n/2}.
\end{equation*}%
Using 3-term recurrence (\ref{_H}) for $q\allowbreak =\allowbreak 1$ and
remembering that $\left[ n\right] _{1}\allowbreak =\allowbreak 1$ we see
that polynomials $\left\{ H_{n}\left( \frac{x-\rho y}{\sqrt{1-\rho ^{2}}}%
\right) \right\} $ satisfy the following 3-term recurrence: 
\begin{equation*}
H_{n+1}\left( \frac{x-\rho y}{\sqrt{1-\rho ^{2}}}\right) \allowbreak
=\allowbreak \frac{x-\rho y}{\sqrt{1-\rho ^{2}}}H_{n}\left( \frac{x-\rho y}{%
\sqrt{1-\rho ^{2}}}\right) -nH_{n-1}\left( \frac{x-\rho y}{\sqrt{1-\rho ^{2}}%
}\right) ,
\end{equation*}%
which after multiplying both sides $\left( 1-\rho ^{2}\right) ^{(n+1)/2}$
and applying definition of $\tilde{P}_{n}$ reduces to 
\begin{equation}
\tilde{P}_{n+1}(x)\allowbreak =\allowbreak \left( x-\rho y\right) \tilde{P}%
_{n}(x)-n(1-\rho ^{2})\tilde{P}_{n-1}(x),  \label{pom}
\end{equation}%
with initial conditions $\tilde{P}_{-1}\left( x\right) \allowbreak
=\allowbreak 0,$ $\tilde{P}_{0}\left( x\right) \allowbreak =\allowbreak 1.$
Now let us consider equation (\ref{ASC}) for $q\allowbreak =\allowbreak 1$
(it is a polynomial in $q)$. It is elementary to notice that $P_{n}(x|y,\rho
,1)$ satisfies 3-term recurrence (\ref{pom}) with the same initial
conditions.

ii) Follows formula (\ref{ASC_HB}) presented below and the properties of
polynomials $B_{n}$ also listed below. In particular the fact that $%
B_{n}(x|0)\allowbreak =\allowbreak 0$ for $n\geq 3.$
\end{proof}

This nice interpretation presented in the Remark \ref{special} is the second
reason why we prefer to work with the polynomials $P_{n}$.

It is known (again e.g. \cite{bms}, \cite{IA} or \cite{Szab5}) that 
\begin{equation}
\sum_{i=0}^{\infty }\frac{t^{i}}{\left[ i\right] _{q}!}P_{i}\left( x|y,\rho
,q\right) \allowbreak =\allowbreak \frac{\varphi \left( x|t,q\right) }{%
\varphi \left( y|\rho t,q\right) },  \label{chASC}
\end{equation}%
where function $\varphi \left( x|t,q\right) $ is given by (\ref{chH}) and
that 
\begin{equation}
\int_{S\left( q\right) }P_{n}\left( x|y,\rho ,q\right) P_{m}\left( x|y,\rho
,q\right) f_{CN}\left( x|y,\rho ,q\right) dx=\left( \rho ^{2}\right) _{n} 
\left[ n\right] _{q}!\delta _{m,n},  \label{pnpm}
\end{equation}%
where 
\begin{equation}
f_{CN}\left( x|y,\rho ,q\right) \allowbreak =\allowbreak f_{N}\left(
x|q\right) \frac{\left( \rho ^{2}\right) _{\infty }}{\prod_{k=0}^{\infty
}w\left( x,y,\rho q^{k}|q\right) },  \label{fCN}
\end{equation}%
with $w\left( x,y,\rho |q\right) \ $defined by: 
\begin{equation}
w\left( x,y,\rho |q\right) =(1-\rho ^{2})^{2}-(1-q)\rho (1+\rho
^{2})xy+(1-q)\rho ^{2}(x^{2}+y^{2}).  \label{w_k}
\end{equation}

In (\ref{chASC}) and (\ref{fCN}) we assume $\left\vert q\right\vert
<1,\allowbreak x,y\allowbreak \in \allowbreak S\left( q\right) ,\left\vert
\rho \right\vert \leq 1,\left\vert t\sqrt{1-q}\right\vert <1$ assuring that
the series and infinite product are absolutely convergent.

We have also the following relationship between the Poisson kernel built of $%
q-$Hermite polynomials and the density $f_{CN}:$ 
\begin{equation}
\frac{\left( \rho ^{2}\right) _{\infty }}{\prod_{k=0}^{\infty }w\left(
x,y,\rho q^{k}|q\right) }=\sum_{i\geq 1}\frac{\rho ^{i}}{\left[ i\right]
_{q}!}H_{i}\left( x|q\right) H_{i}\left( y|q\right) .  \label{mehler}
\end{equation}%
There are many proofs of (\ref{mehler}). For some see e.g. \cite{IA}, \cite%
{bressoud}, \cite{Expansion}. Convergence here is absolute for $\left\vert
\rho \right\vert <1,$ $x,y\in S\left( q\right) .$ (\ref{mehler}) is known
under the name the Poisson--Mehler expansion formula. Generally discussion
of the convergence conditions for expressions involving $q-$Hermite
polynomials is done in \cite{Expansion} or \cite{AW}.

Again as before we have the following remark where special values of
parameter $q$ are considered.

\begin{remark}
\label{fCN_specal}i) For $q\allowbreak =\allowbreak 1$ we have $f_{CN}\left(
x|y,\rho ,1\right) \allowbreak =\allowbreak \exp \left( -\frac{(x-\rho y)^{2}%
}{2(1-\rho ^{2})}\right) /\sqrt{2\pi (1-\rho ^{2})}.$

ii) For $q\allowbreak =\allowbreak 0$ we have $f_{CN}\left( x|y,\rho
,0\right) \allowbreak =\allowbreak \frac{(1-\rho ^{2})\sqrt{4-x^{2}}}{2\pi
w(x,y,\rho |0)}$ (so called Kesten--McKay density).
\end{remark}

\begin{proof}
i) Rigorous proof that $f_{CN}\left( x|y,\rho ,q\right) \longrightarrow \exp
\left( -\frac{(x-\rho y)^{2}}{2(1-\rho ^{2})}\right) /\sqrt{2\pi (1-\rho
^{2})}$ as $q\longrightarrow 1^{-}$ can be found in \cite{ISV87}. Less
formal and more intuitive argument is the following. By Remark \ref{special}
i) we deduce that all moments tend (as $q\longrightarrow 1^{-})$ to the
moments of the limiting distribution, hence distribution with the density $%
f_{CN}$ must tend weakly to a distribution with the density $\exp \left( -%
\frac{(x-\rho y)^{2}}{2(1-\rho ^{2})}\right) /\sqrt{2\pi (1-\rho ^{2})}.$

ii) We simply set $q=0$ in (\ref{fCN}).
\end{proof}

Furthermore we will need some auxiliary polynomials that are in fact related
to $H_{n}\left( x|q^{-1}\right) ,$ $-1<q\leq 1.$ Namely let us consider
polynomials $\left\{ B_{n}\right\} _{n\geq -1}$ satisfying the following
3-term recurrence%
\begin{equation*}
B_{n+1}\left( x|q\right) =-q^{n}xB_{n}\left( x|q\right)
+q^{n-1}[n]_{q}B_{n-1}\left( x|q\right) :\ (n\geq 0)
\end{equation*}%
with the usual initial conditions $B_{-1}=0,B_{0}=1$. To support intuition,
let us remark (following \cite{bms}), that 
\begin{eqnarray*}
B_{n}\left( x|1\right) \allowbreak &=&\allowbreak i^{n}H_{n}\left( ix\right)
, \\
\text{and }B_{n}\left( x|0\right) \allowbreak &=&\allowbreak 0,
\end{eqnarray*}%
for $n\geq 3,$ with $B_{0}(x|0)=1,$ $B_{1}\left( x|0\right) \allowbreak
=\allowbreak -x,$ $B_{2}\left( x|0\right) \allowbreak =\allowbreak 1.$

It was shown in \cite{bms} that:%
\begin{equation}
\forall n\geq 0:P_{n}(x|y,\rho ,q)=\sum_{i=0}^{n}\QATOPD[ ] {n}{i}_{q}\rho
^{n-i}B_{n-i}\left( y|q\right) H_{i}\left( x|q\right) .  \label{ASC_HB}
\end{equation}

We will use this formula to extend definition of ASC polynomials for $%
\left\vert \rho \right\vert \geq 1.$ Of course for $\left\vert \rho
\right\vert \geq 1$ ASC polynomials are not orthogonal with respect to a
positive measure.

Now let us recall the idea of expansion of the ratio of two densities in a
series of polynomials orthogonal with respect to the one of them presented
with many examples in \cite{Expansion}. More precisely let us assume that we
have two positive probability measures $d\alpha $ and $d\beta $ with
densities respectively $A\left( x\right) $ and $B\left( x\right) $. Moreover
let us assume that $\left\{ a_{n}\left( x\right) \right\} _{n\geq 0}$ and $%
\left\{ b_{n}\left( x\right) \right\} _{n\geq 0}$ are sets of monic
polynomials orthogonal with respect to $d\alpha $ and $d\beta $
respectively. Besides assume that $\limfunc{supp}\beta \allowbreak
=\allowbreak \limfunc{supp}\alpha $.

Further let us also assume that we know the so called 'connection
coefficients' between set $\left\{ b_{n}\right\} $ and $\left\{
a_{n}\right\} $. More precisely let us assume that we know the numbers $%
\gamma _{k,n}$ such that for every $n\geq 0:$%
\begin{equation*}
a_{n}\left( x\right) \allowbreak =\allowbreak \sum_{k=0}^{n}\gamma
_{k,n}b_{k}\left( x\right) 
\end{equation*}%
for every $x\in \mathbb{R}$. Besides assume that $\int_{\limfunc{supp}\left(
B\right) }(B^{2}\left( x\right) /A^{2}\left( x\right) )d\alpha (x)<\infty ,$
then we have 
\begin{equation}
B\left( x\right) \allowbreak =\allowbreak A\left( x\right)
\sum_{n=0}^{\infty }\frac{\gamma _{0,n}\hat{b}_{0}}{\hat{a}_{n}}a_{n}\left(
x\right) ,  \label{series}
\end{equation}%
where $\hat{a}_{n}\allowbreak =\allowbreak \int_{\limfunc{supp}\left(
A\right) }a_{n}\left( x\right) ^{2}d\alpha \left( x\right) ,$ similarly for
the polynomials $b_{n}.$ Convergence in (\ref{series}) is $L^{2}$
convergence, however if coefficients $\frac{\gamma _{0,n}\hat{b}_{0}}{\hat{a}%
_{n}}$ are such that $\sum_{n\geq 0}\left( \frac{\gamma _{0,n}\hat{b}_{0}}{%
\hat{a}_{n}}\right) ^{2}\log ^{2}n$\allowbreak $<\allowbreak \infty ,$ then
by the Rademacher-Menshov theorem we have almost pointwise, absolute
convergence. In most cases interesting in the $q-\func{series}$ theory this
condition is trivially satisfied hence all expansions we are going to
consider will be absolutely almost pointwise convergent.

Hence finding connection coefficients between two sets of orthogonal
polynomials becomes crucial.

Returning to the big $q$-Hermite and ASC polynomials we have the following
easy lemma.

\begin{lemma}
\label{connection}i) $\forall n\geq 0,x,y,b,a,q\in \mathbb{R}:$%
\begin{equation*}
H_{n}(x|a,q)\allowbreak =\allowbreak \sum_{j=0}^{n}\QATOPD[ ] {n}{j}%
_{q}P_{j}\left( x|y,\frac{a}{b},q\right) \left( \frac{a}{b}\right)
^{n-j}H_{n-j}\left( y|b,q\right) .
\end{equation*}

ii) $\forall n\geq 0,x,y,\rho ,a,q\in \mathbb{R}:$%
\begin{equation*}
P_{n}\left( x|y,\rho ,q\right) \allowbreak =\allowbreak \sum_{i=0}^{n}\QATOPD%
[ ] {n}{i}_{q}\rho ^{n-i}B_{n-i}\left( y|a/\rho ,q\right) H_{i}(x|a,q),
\end{equation*}%
where we denoted $B_{m}\left( x|b,q\right) \overset{df}{=}\sum_{j=0}^{m}%
\QATOPD[ ] {m}{j}_{q}b^{m-j}B_{j}\left( x|q\right) .$

iii) $\forall n\geq 0,x,y,\rho _{1},\rho _{2},q\in \mathbb{R}:$%
\begin{equation*}
P_{n}\left( x|y,\rho _{1}\rho _{2},q\right) \allowbreak =\allowbreak
\sum_{i=0}^{n}\QATOPD[ ] {n}{i}_{q}\rho _{1}^{n-i}P_{i}\left( x|z,\rho
_{1},q\right) P_{n-i}\left( z|y,\rho _{2},q\right) .
\end{equation*}

iv) $\forall n\geq 1,\left\vert q\right\vert <1:$%
\begin{eqnarray*}
\max_{x\in S\left( q\right) }\left\vert H_{n}\left( x|q\right) \right\vert
\allowbreak &\leq &\left( 1-q\right) ^{-n/2}r_{n}\left( 1|q\right) , \\
\max_{x\in S\left( q\right) }\left\vert H_{n}(x|a,q)\right\vert \allowbreak
&\leq &\allowbreak \left( -\left\vert a\right\vert \sqrt{1-q}\right)
_{\infty }\left( 1-q\right) ^{-n/2}r_{n}\left( 1|q\right) ,
\end{eqnarray*}%
$\,$

where $\allowbreak r_{n}\left( x|q\right) $ is given by $r_{n}\left(
x|q\right) \allowbreak =\allowbreak \sum_{i=0}^{n}\QATOPD[ ] {n}{i}%
_{q}x^{i}. $

v) $\forall \left\vert q\right\vert <1,x,y\in S\left( q\right) :0<\frac{%
\left( \rho ^{2}\right) _{\infty }}{\left( -\left\vert \rho \right\vert
\right) _{\infty }^{4}}\footnote{%
This form of the lower bound was sugested by the referee}\leq \frac{%
f_{CN}\left( x|y,\rho ,q\right) }{f_{N}\left( x|q\right) }\leq \frac{\left(
\rho ^{2}\right) _{\infty }}{\left( \left\vert \rho \right\vert \right)
_{\infty }^{4}}$
\end{lemma}

Easy and not very interesting proof of this lemma has been moved to Section %
\ref{dowody}.

\section{Main results}

Now using the described above idea of density expansion and Lemma \ref%
{connection} we have the following theorem.

\begin{theorem}
\label{glowne} For $\left\vert b\right\vert >\left\vert a\right\vert
,\left\vert q\right\vert <1,$ $x,y\in S\left( q\right) $ we have

i)%
\begin{equation}
0\leq \sum_{n\geq 0}\frac{a^{n}}{\left[ n\right] _{q}!b^{n}}%
H_{n}(x|a,q)H_{n}\left( y|b,q\right) \allowbreak =\allowbreak \left(
a^{2}/b^{2}\right) _{\infty }\prod_{k=0}^{\infty }\frac{\left(
1-(1-q)xaq^{k}+(1-q)a^{2}q^{2k}\right) }{w\left( x,y,q^{k}a/b|q\right) },
\label{jadro}
\end{equation}

ii) 
\begin{equation*}
1/\sum_{n\geq 0}\frac{a^{n}}{\left[ n\right] _{q}!b^{n}}H_{n}(x|a,q)H_{n}%
\left( y|b,q\right) =\sum_{n\geq 0}\frac{a^{n}}{\left[ n\right]
_{q}!b^{n}\left( a^{2}/b^{2}\right) _{n}}B_{n}\left( y|b,q\right)
P_{n}\left( x|y,a/b,q\right) .
\end{equation*}

iii) For $\left\vert \rho _{1}\right\vert ,\left\vert \rho _{2}\right\vert
,\left\vert q\right\vert <1,$ $x,y,z\in S\left( q\right) $ 
\begin{equation}
0\leq \sum_{n\geq 0}\frac{\rho _{1}^{n}}{\left[ n\right] _{q}!\left( \rho
_{1}^{2}\rho _{2}^{2}\right) _{n}}P_{n}\left( x|y,\rho _{1}\rho
_{2},q\right) P_{n}\left( z|y,\rho _{2},q\right) =\frac{\left( \rho
_{1}^{2}\right) _{\infty }}{\left( \rho _{1}^{2}\rho _{2}^{2}\right)
_{\infty }}\prod_{k=0}^{\infty }\frac{w\left( x,y,q^{k}\rho _{1}\rho
_{2}|q\right) }{w\left( x,z,q^{k}\rho _{1}|q\right) }.  \label{kerASC}
\end{equation}
\end{theorem}

\begin{proof}
i) Let us apply assertion of the Lemma \ref{connection} with $\rho
\allowbreak =\allowbreak a/b.$ We get then: 
\begin{equation*}
H_{n}(x|a,q)\allowbreak =\allowbreak \sum_{i=0}^{n}\QATOPD[ ] {n}{i}%
_{q}P_{i}\left( x|y,a/b,q\right) \left( a/b\right) ^{n-i}H_{n-i}\left(
y|b,q\right) .
\end{equation*}
Hence coefficient $\gamma _{0,n}\allowbreak =\allowbreak \left( a/b\right)
^{n}H_{n}\left( y|b,q\right) .$ Thus applying described above idea of
expanding ratio of two densities we get%
\begin{equation}
f_{CN}\left( x|y,a/b,q\right) \allowbreak =\allowbreak f_{N}\left(
x|q\right) \varphi (x|a,q)\sum_{n\geq 0}\frac{\left( a/b\right) ^{n}}{\left[
n\right] _{q}!}H_{n}(x|a,q)H_{n}\left( y|b,q\right) .  \label{exp1}
\end{equation}%
Since we have $\int_{S\left( q\right) }H_{n}^{2}(x|a,q)f_{N}\left(
x|q\right) \varphi (x|a,q)dx\allowbreak =\allowbreak \left[ n\right] _{q}!$
we apply (\ref{fCN}).

ii) To get ii) we reason in the similar way this time using $P_{n}\left(
x|y,a/b,q\right) \allowbreak =\allowbreak \sum_{i=0}^{n}\QATOPD[ ] {n}{i}%
_{q}(a/b)^{n-i}B_{n-i}\left( y|b,q\right) H_{i}(x|a,q)$ and keeping in mind (%
\ref{pnpm}). As far as convergence of both series is concerned we see that
for $\left\vert \rho \right\vert ,\left\vert q\right\vert <1$ and $x,y$ $\in
S_{q}$ function $g(x|y,\rho ,q)\allowbreak =\allowbreak f_{CN}\left(
x|y,\rho ,q\right) /f_{bN}\left( x|q\right) \allowbreak $ is both bounded
and 'cut away from zero' (compare (\ref{fbN}) and (\ref{fCN})). Hence its
square as well as reciprocal of this square are integrable on compact
interval $S_{q}.$

iii) We reason in the similar way using assertion iii) of the Lemma \ref%
{connection}. This time coefficient $\gamma _{0,n}\allowbreak =\allowbreak
\rho _{1}^{n}P_{n}\left( z|y,\rho _{2},q\right) .$ Since we have (\ref{pnpm}%
) we can write 
\begin{equation}
f_{CN}\left( x|z,\rho _{1},q\right) \allowbreak =\allowbreak f_{CN}\left(
x|y,\rho _{1}\rho _{2},q\right) \sum_{n\geq 0}\frac{\rho _{1}^{n}}{\left[ n%
\right] _{q}!\left( \rho _{1}^{2}\rho _{2}^{2}\right) _{n}}P_{n}\left(
x|y,\rho _{1}\rho _{2},q\right) P_{n}\left( z|y,\rho _{2},q\right) .
\label{exp2}
\end{equation}%
In all three cases convergence is absolute and almost sure since we have
assertion v) of Lemma \ref{connection}.
\end{proof}

Notice that the right hand side of (\ref{kerASC}) can be written as $%
f_{CN}(x|z,\rho _{1},q)/f_{CN}\left( x|y,\rho _{1}\rho _{2},q\right) $ which
would enable to pass easily with $q$ to $1^{-}.$ Keeping this in mind we
have the following statement as a corollary:

\begin{corollary}
i) Setting $q\allowbreak =\allowbreak 1$ we get for $x,y\in \mathbb{R}$ and $%
\left\vert a\right\vert <\left\vert b\right\vert $:%
\begin{equation*}
\sum_{n\geq 0}\frac{\rho ^{n}}{n!}H_{n}\left( x-a\right) H_{n}\left(
y-b\right) \allowbreak =\allowbreak \frac{1}{\sqrt{1-\rho ^{2}}}\exp \left( -%
\frac{\left( x-\rho y\right) ^{2}}{2(1-\rho ^{2})}+\frac{(x-a)^{2}}{2}%
\right) .
\end{equation*}%
$\allowbreak $

ii) Setting $q\allowbreak =\allowbreak 0$ we get for $x,y\in \lbrack -2,2]$ $%
,$ $1>\left\vert b\right\vert >\left\vert a\right\vert $ 
\begin{eqnarray*}
&&\sum_{n\geq 0}^{n}\rho ^{n}\left( U_{n}\left( x/2\right) -aU_{n-1}\left(
x/2\right) \right) \left( U_{n}\left( y/2\right) -bU_{n-1}\left( y/2\right)
\right) \\
&=&\left( 1-\rho ^{2}\right) \frac{\left( 1-ax+a^{2}\right) }{w(x,y,\rho |0)}%
.
\end{eqnarray*}

In both formulae above we denoted for simplicity $\rho \allowbreak
=\allowbreak a/b.$

iii) Letting $q\longrightarrow 1^{-}$ in (\ref{kerASC}) we get for $x,y,z\in 
\mathbb{R}$ and $\left\vert \rho _{1}\right\vert ,\left\vert \rho
_{2}\right\vert <1$:%
\begin{eqnarray*}
&&\sum_{n=0}^{\infty }\frac{\rho _{1}^{n}\left( 1-\rho _{2}^{2}\right) ^{n/2}%
}{n!\left( 1-\rho _{1}^{2}\rho _{2}^{2}\right) ^{n/2}}H_{n}\left( \frac{%
x-\rho _{1}\rho _{2}y}{\sqrt{(1-\rho _{1}^{2}\rho _{2}^{2})}}\right)
H_{n}\left( \frac{z-\rho _{2}y}{\sqrt{(1-\rho _{2}^{2})}}\right)  \\
&=&\sqrt{\frac{1-\rho _{1}^{2}\rho _{2}^{2}}{1-\rho _{1}^{2}}}\exp \left( -%
\frac{(x-\rho _{1}z)^{2}}{2(1-\rho _{1}^{2})}+\frac{(x-\rho _{1}\rho
_{2}y)^{2}}{2(1-\rho _{1}^{2}\rho _{2}^{2})}\right) .
\end{eqnarray*}

iv) Setting $q\allowbreak =\allowbreak 0$ in (\ref{kerASC}) we get for $%
x,y,z\in \lbrack -2,2]$ and $\left\vert \rho _{1}\right\vert ,\left\vert
\rho _{2}\right\vert <1:$%
\begin{eqnarray*}
&&1+\frac{1}{\left( 1-\rho _{1}^{2}\rho _{2}^{2}\right) }\sum_{n\geq 1}\rho
_{1}^{n}P_{n}\left( x|y,\rho _{1}\rho _{2},0\right) P_{n}\left( z|y,\rho
_{2},0\right) \\
&=&\frac{(1-\rho _{1}^{2})w\left( x,y,\rho _{1}\rho _{2}|0\right) }{(1-\rho
_{1}^{2}\rho _{2}^{2})w(x,z,\rho _{1}|0)}.
\end{eqnarray*}
\end{corollary}

\begin{proof}
i) First we consider (\ref{ASC}) for $q\allowbreak =\allowbreak 1$ and
notice that $\varphi \left( x|a,1\right) \allowbreak =\allowbreak \exp
(ax-a^{2}/2),$ $f_{N}(x)\allowbreak =\allowbreak \exp (-x^{2}/2)/\sqrt{2\pi }%
.$ By Lemma \ref{connection}, i) and Remark \ref{bH_special}, i) we have:%
\begin{equation*}
H_{n}(x-a)\allowbreak =\allowbreak \sum_{j=0}^{n}\binom{n}{j}P_{j}(x|y,\frac{%
a}{b},1)\left( \frac{a}{b}\right) ^{n-j}H_{n-j}\left( y-b\right) .
\end{equation*}%
Now we apply (\ref{series}) which proves (\ref{exp1}) for $q\allowbreak
=\allowbreak 1$ and consequently implies assertion. ii) We use (\ref{jadro})
and Remark \ref{bH_special}. iii) First we use Lemma \ref{connection}, iii)
with $q\allowbreak =\allowbreak 1$ and Remark \ref{special}, i). Then we use
(\ref{series}) and show (\ref{exp2}) which implies the assertion. iv) Again
we use (\ref{kerASC}) to get left hand side while (\ref{exp2}) and the
second assertion of Remark \ref{fCN_specal} to get the right hand side.
\end{proof}

\begin{remark}
In \cite{suslov96} the following formula (14.14) (expressed in terms of
polynomials $h_{n}(x|a,q)$ with slightly different definition of big $q-$%
Hermite polynomials) was given:%
\begin{gather*}
\sum_{n\geq 0}\left( t\frac{a}{b}\right) ^{n}\frac{1}{\left( q\right) _{n}}%
h_{n}(x|a,q)h_{n}\left( y|b,q\right) \allowbreak =\allowbreak \frac{\left(
a^{2}t^{2}/b^{2},tae^{i\theta },tae^{-i\theta }\right) _{\infty }}{\left(
ate^{i\left( \theta +\phi \right) }/b,tae^{i(\theta -\phi )}/b,tae^{i(\phi
-\theta )}/b,tae^{-i(\theta +\phi )}/b\right) _{\infty }}\times  \\
_{3}\phi _{2}\left( 
\begin{array}{ccc}
t & ate^{i\left( \theta +\phi \right) }/b & ate^{i\left( \theta -\phi
\right) }/b \\ 
a^{2}t^{2}/b^{2} & ate^{i\theta } & 
\end{array}%
;q,ae^{-i\theta }\right) ,
\end{gather*}%
with $x\allowbreak =\allowbreak \cos \theta $ and $y\allowbreak =\allowbreak
\cos \phi $ convergent for $\left\vert x\right\vert ,\left\vert y\right\vert
,\left\vert t\right\vert \leq 1$ and $\left\vert q\right\vert <1.$ If $%
\left\vert b\right\vert <\left\vert a\right\vert ,$ $t\allowbreak
=\allowbreak 1$ we get $\left. _{3}\phi _{2}\left( 
\begin{array}{ccc}
t & ate^{i\left( \theta +\phi \right) }/b & ate^{i\left( \theta -\phi
\right) }/b \\ 
a^{2}t^{2}/b^{2} & ate^{i\theta } & 
\end{array}%
;q,ae^{-i\theta }\right) \right\vert _{t=1}\allowbreak =\allowbreak 1,$
hence when polynomials $h_{n}$ are replaced by polynomials $H_{n},$ $a,$ $b,$
$x$ and $y$ by $a\sqrt{1-q},$ $b\sqrt{1-q},$ $x\sqrt{1-q}/2$ and $y\sqrt{1-q}%
/2$ respectively, we get left hand side of (\ref{jadro}). To get right hand
side of (\ref{jadro}) we must use the fact that $\cos \theta $ and $\cos
\phi $ should be replaced $x\sqrt{1-q}/2$ and $y\sqrt{1-q}/2$ respectively
and also the following observations:%
\begin{gather*}
(1-ae^{i\theta }\sqrt{1-q}q^{k})(1-ae^{-i\theta }\sqrt{1-q}q^{k})\allowbreak
=\allowbreak 1-\sqrt{1-q}aq^{k}\cos \theta +(1-q)a^{2}, \\
(1-\frac{a}{b}e^{i\left( \theta +\phi \right) }q^{k})(1-\frac{a}{b}%
e^{i\left( \theta -\phi \right) }q^{k})(1-\frac{a}{b}e^{i\left( -\theta
+\phi \right) }q^{k})(1-\frac{a}{b}e^{-i\left( \theta +\phi \right)
}q^{k})=w\left( x,y,q^{k}\frac{a}{b}|q\right) .
\end{gather*}
However our proof of (\ref{jadro}) is much simpler than the proof presented
in \cite{suslov96}. Besides we obtain asymmetric kernel expansion of the
reciprocal of (\ref{jadro}).
\end{remark}

\begin{remark}
In \cite{suslov96} there are formulae 14.5, and 14.8 (expressed in terms of
polynomials $p_{n}\left( x|a,b,q\right) $ compare (\ref{pn})) for summing
kernels of the form \newline
$\sum_{n\geq 0}\frac{\left( ab\right) _{n}}{\left( q\right) _{n}}\left( t%
\frac{\alpha }{a}\right) ^{n}p_{n}\left( x|a,b,q\right) p_{n}\left( y|\alpha
,\beta ,q\right) $ and $\sum_{n\geq 0}\frac{\left( q\right) _{n}\left(
ta\alpha \right) ^{n}}{\left( ab\right) _{n}}p_{n}\left( x|a,b,q\right)
p_{n}\left( y,\alpha ,\beta ,q\right) $ under condition that $\alpha \beta
=ab$ and of course $\left\vert t\right\vert \leq 1$. It was expressed in
terms of the basic hypergeometric function $_{8}\phi _{7}.$ Positivity of
those kernels was shown only in some particular special cases. Slightly more
general form of kernels involving ASC polynomials are presented in \cite%
{IsStan02} .

Our kernel (\ref{kerASC}) is different. We do not need assumption $\alpha
\beta \allowbreak =\allowbreak ab$ which, as indicated above (compare (\ref%
{pn}) and (\ref{Pn})), leads to the condition : $\rho _{1}^{2}\allowbreak
\rho _{2}^{2}=\allowbreak \rho _{2}^{2}$. Recall that in our setting
parameters $a,b$ are related to parameters $\rho $, $y$ by $ab\allowbreak
=\allowbreak \rho ^{2}$ and $a+b\allowbreak =\allowbreak \rho y\sqrt{1-q}.$
We also assume that parameters $a,b$, $\alpha ,\beta $ are related to one
another however in a different way namely $\left( a+b\right) /\sqrt{ab}%
\allowbreak =\allowbreak \left( \alpha +\beta \right) /\sqrt{\alpha \beta }.$
Thus our result although connected to known results is different and new, it
was obtained by much simpler argument. Besides we have positivity of our
kernels again 'for free'.
\end{remark}

\begin{remark}
Notice that putting $\rho _{2}\allowbreak =\allowbreak 0$ in (\ref{kerASC})
we get the Poisson--Mehler formula (i.e. formula (\ref{mehler})) since $%
P_{n}\left( x|y,0,q\right) \allowbreak =\allowbreak H_{n}\left( x|q\right) $%
. Hence Theorem \ref{glowne} provides yet another (and very simple,
elementary) proof of the Poisson--Mehler formula. Similarly passing with $a$
and $b$ to zero in such a way that $a/b\longrightarrow \rho $ we see that (%
\ref{jadro}) is another generalization of the Poisson--Mehler formula.
\end{remark}

\begin{remark}
If one formally extends the definition of ASC polynomials by using assertion
ii) of Lemma \ref{connection} say with $a\allowbreak =\allowbreak 0$ for $%
\left\vert \rho \right\vert \geq 1$ (thus loosing the fact that they are
orthogonal with respect to a positive measure) then assertion iii) of
Theorem \ref{glowne} can be rewritten in a more symmetric way (after
redefining $\rho _{1}$ and $\rho _{2})$ in the following form: for $%
\left\vert \rho _{1}\right\vert ,\left\vert \rho _{2}\right\vert ,\left\vert
q\right\vert <1,$ $x,y,z\in S\left( q\right) $ 
\begin{equation}
0\leq \sum_{n\geq 0}\frac{\rho _{1}^{n}}{\left[ n\right] _{q}!\left( \rho
_{2}^{2}\right) _{n}}P_{n}\left( x|y,\rho _{2},q\right) P_{n}\left( z|y,%
\frac{\rho _{2}}{\rho _{1}},q\right) =\frac{\left( \rho _{1}^{2}\right)
_{\infty }}{\left( \rho _{2}^{2}\right) _{\infty }}\prod_{k=0}^{\infty }%
\frac{w\left( x,z,\rho _{2}q^{k}|q\right) }{w\left( x,y,\rho
_{1}q^{k}|q\right) },  \label{genASC}
\end{equation}%
from which directly follows the following inversion of the kernel formula:
for $\left\vert \rho _{1}\right\vert ,\left\vert \rho _{2}\right\vert
,\left\vert q\right\vert <1,$ $x,y,z\in S\left( q\right) :$%
\begin{equation}
1=\sum_{n\geq 0}\frac{\rho _{1}^{n}}{\left[ n\right] _{q}!\left( \rho
_{2}^{2}\right) _{n}}P_{n}\left( x|y,\rho _{2},q\right) P_{n}\left( z|y,%
\frac{\rho _{2}}{\rho _{1}},q\right) \times \sum_{n\geq 0}\frac{\rho _{2}^{n}%
}{\left[ n\right] _{q}!\left( \rho _{1}^{2}\right) _{n}}P_{n}\left( x|z,\rho
_{1},q\right) P_{n}\left( y|z,\frac{\rho _{1}}{\rho _{2}},q\right) .
\label{nice}
\end{equation}
\end{remark}

\begin{remark}
Applying just for fun of checking the idea of expansion of the ratio of
densities to the formulae (\ref{bigH}) and (\ref{bigH2}) treated as the
'connection coefficient' formulae we get respectively: 
\begin{eqnarray}
\frac{1}{\varphi (x|a,q)}\allowbreak  &=&\allowbreak \sum_{n\geq 0}\frac{%
\left( -a\right) ^{n}}{\left[ n\right] _{q}!}q^{\binom{n}{2}}H_{n}(x|a,q),
\label{_1} \\
\varphi (x|a,q)\allowbreak  &=&\allowbreak \sum_{n\geq 0}\frac{a^{n}}{\left[
n\right] _{q}!}H_{n}\left( x|q\right) .  \label{_2}
\end{eqnarray}%
Thus having yet another proof of the formula defining generating function of 
$q-$Hermite polynomials (\ref{_2}). As far as the formula (\ref{_1}) is
concerned in \cite{Koek} we find formula (3.18.14) which gives the sum $%
\sum_{n\geq 0}\frac{\left( -t\right) ^{n}}{\left( q\right) _{n}!}q^{\binom{n%
}{2}}h_{n}(x|a,q)\allowbreak $ in terms of hypergeometric function $_{1}\phi
_{1}.$ Thus (\ref{_1}) gives special value of this generating function for $%
t\allowbreak =\allowbreak a.$
\end{remark}

\section{Proofs\label{dowody}}

\begin{proof}[Proof of Lemma \protect\ref{connection}]
We shall prove our identities for $\left\vert q\right\vert ,\left\vert \rho
\right\vert <1.$ Then since both the left hand side and the right hand side
of the above mentioned identity are polynomials in $q$ and $\rho $ of order
at most $n$ we will deduce that the identity is true for all $q$ and $\rho .$
The proofs below are so to say by a 'direct method'.

However they can also be derived by 'characteristic function method' if one
noticed the following: 
\begin{eqnarray*}
\varphi _{bH}(x|t,a,q)\allowbreak  &=&\sum_{j\geq 0}\frac{t^{n}}{\left[ n%
\right] _{q}!}H_{n}(x|a,q)=\varphi \left( x|t,q\right) \left( (1-q)at\right)
_{\infty }, \\
\varphi _{B}\left( x|t,q\right)  &=&\sum_{j\geq 0}\frac{t^{n}}{\left[ n%
\right] _{q}!}B_{n}(x|q)\allowbreak =\allowbreak \frac{1}{\varphi \left(
x|t,q\right) },~\varphi _{P}\left( x|t,y,\rho ,q\right) \allowbreak
=\allowbreak \frac{\varphi (x|t,q)}{\varphi \left( y|\rho t,q\right) } \\
\varphi _{bB}\left( x|t,a,q\right)  &=&\sum_{j\geq 0}\frac{t^{n}}{\left[ n%
\right] _{q}!}B_{n}(x|a,q)\allowbreak =\allowbreak \frac{1}{\varphi \left(
x|t,q\right) \left( (1-q)at\right) _{\infty }}
\end{eqnarray*}%
which are modified versions of the characteristic functions given in the
literature: $\varphi _{bH}(x|t,a,q)$ is a modification of appropriate
formula from \cite{Koek} and $\varphi _{B}\left( x|t,q\right) $ from \cite%
{bms}. Then for example assertion i) follows identity 
\begin{equation*}
\varphi _{P}\left( x|t,y,\frac{a}{b},q\right) \varphi _{bH}\left( x|\frac{a}{%
b}t,a,q\right) \allowbreak =\allowbreak \varphi _{bH}\left( x|t,a,q\right) .
\end{equation*}%
Similarly for the other assertions.

i) Recall that in \cite{IRS99} formula (4.7) it is a 'connection
coefficient' formula between polynomials $\left\{ h_{n}\right\} $ and $%
\left\{ p_{n}\right\} ,$ which can be easily rewritten in terms of $q-$%
Hermite and ASC polynomials using (\ref{Pn}). Formula (4.7) of \cite{IRS99}
thus now reads: 
\begin{equation}
H_{n}\left( x|q\right) \allowbreak =\allowbreak \sum_{i=0}^{n}\QATOPD[ ] {n}{%
i}_{q}\rho ^{n-i}H_{n-i}\left( y|q\right) P_{i}\left( x|y,\rho ,q\right) .
\label{HnaP}
\end{equation}%
Let us use also (\ref{bigH}). Then we get: 
\begin{eqnarray*}
H_{n}(x|a,q) &=&\sum_{i=0}^{n}\QATOPD[ ] {n}{i}_{q}(-1)^{i}q^{\binom{i}{2}%
}a^{i}H_{n-i}\left( x|q\right) \\
&=&\sum_{i=0}^{n}\QATOPD[ ] {n}{i}_{q}(-1)^{i}q^{\binom{i}{2}%
}a^{i}\sum_{j=0}^{n-i}\QATOPD[ ] {n-i}{j}_{q}\rho ^{n-i-j}H_{n-i-j}\left(
y|q\right) P_{j}\left( x|y,\rho ,q\right) \\
&=&\sum_{j=0}^{n}\QATOPD[ ] {n}{j}_{q}P_{j}\left( x|y,\rho ,q\right) \rho
^{n-j}\sum_{i=0}^{n-j}\QATOPD[ ] {n-j}{i}(-1)^{i}q^{\binom{i}{2}}\left( 
\frac{a}{\rho }\right) ^{i}H_{n-i-j}\left( y|q\right) \\
&=&\sum_{j=0}^{n}\QATOPD[ ] {n}{j}_{q}P_{j}\left( x|y,\rho ,q\right) \rho
^{n-j}H_{n-j}\left( y;a/\rho |q\right) .
\end{eqnarray*}%
Now denote $a/\rho =b.$

ii) We will use (\ref{ASC_HB}) and (\ref{bigH2})$.$ We have%
\begin{eqnarray*}
P_{n}\left( x|y,\rho ,q\right) \allowbreak &=&\allowbreak \sum_{i=0}^{n}%
\QATOPD[ ] {n}{i}_{q}\rho ^{n-i}B_{n-i}\left( y|q\right) H_{i}\left(
x|q\right) \allowbreak \allowbreak \\
&=&\sum_{i=0}^{n}\QATOPD[ ] {n}{i}_{q}\rho ^{n-i}B_{n-i}\left( y|q\right)
\sum_{k=0}^{i}\QATOPD[ ] {i}{k}_{q}a^{i-k}H_{k}(x|a,q) \\
&=&\sum_{k=0}^{n}\QATOPD[ ] {n}{k}_{q}H_{k}(x|a,q)\sum_{i=k}^{n}\QATOPD[ ] {%
n-k}{i-k}_{q}\rho ^{n-i}B_{n-i}\left( y|q\right) a^{i-k} \\
&=&\sum_{k=0}^{n}\QATOPD[ ] {n}{k}_{q}H_{k}(x|a,q)\rho ^{n-k}\sum_{m=0}^{n-k}%
\QATOPD[ ] {n-k}{m}_{q}(a/\rho )^{m}B_{n-k-m}\left( y|q\right)
\end{eqnarray*}

iii) Keeping in mind (\ref{ASC_HB}) and (\ref{HnaP}) we get: 
\begin{gather*}
P_{n}\left( x|y,\rho ,q\right) \allowbreak =\allowbreak \sum_{i=0}^{n}\QATOPD%
[ ] {n}{i}_{q}\rho ^{n-i}B_{n-i}\left( y|q\right) H_{i}\left( x|q\right)
\allowbreak \\
=\sum_{i=0}^{n}\QATOPD[ ] {n}{i}_{q}\rho ^{n-i}B_{n-i}\left( y|q\right)
\sum_{j=0}^{i}\QATOPD[ ] {i}{j}_{q}r^{i-j}H_{i-j}\left( z|q\right)
P_{j}\left( x|z,r,q\right) \\
=\sum_{j=0}^{n}\QATOPD[ ] {n}{j}_{q}P_{j}\left( x|z,r,q\right) \sum_{i=j}^{n}%
\QATOPD[ ] {n-j}{i-j}_{q}\rho ^{n-i}r^{i-j}B_{n-i}\left( y|q\right)
H_{i-j}\left( z|q\right) \allowbreak \\
=\sum_{j=0}^{n}\QATOPD[ ] {n}{j}_{q}P_{j}\left( x|z,r,q\right)
\sum_{s=0}^{n-j}\QATOPD[ ] {n-j}{s}_{q}\rho ^{n-j-s}r^{s}B_{n-j-s}\left(
y|q\right) H_{s}\left( z|q\right) .
\end{gather*}%
Let us extend definition of polynomials $P_{n}\left( x|y,t,q\right) $ for $%
\left\vert t\right\vert \geq 1$ by defining $P_{n}\left( x|y,t,q\right)
=\sum_{k=0}^{n}\QATOPD[ ] {n}{k}_{q}t^{n-k}B_{n-k}\left( y|q\right)
H_{k}\left( x|q\right) .\allowbreak $

Thus we have: 
\begin{equation*}
P_{n}\left( x|y,\rho ,q\right) =\sum_{j=0}^{n}\QATOPD[ ] {n}{j}%
_{q}r^{n-j}P_{j}\left( x|z,r,q\right) P_{n-j}(z|y,\rho /r,q).
\end{equation*}%
Now let us redefine $\rho $ and $r$ by selecting them so that $\left\vert
r\right\vert >\left\vert \rho \right\vert $ and changing notation $\rho
/r\longrightarrow \rho _{1},$ $r\longrightarrow \rho _{2},$ consequently $%
\rho \longrightarrow \rho _{1}\rho _{2}.$

iv) 
\begin{eqnarray*}
\left\vert H_{n}(x|a,q)\right\vert \allowbreak  &\leq &\allowbreak
\sum_{i=0}^{n}\QATOPD[ ] {n}{i}q^{\binom{i}{2}}\allowbreak \left\vert
a\right\vert ^{i}\allowbreak \left\vert H_{n-i}\left( x|q\right) \right\vert
\allowbreak  \\
&\leq &\allowbreak \sum_{i=0}^{n}\QATOPD[ ] {n}{i}\allowbreak q^{\binom{i}{2}%
}\allowbreak \left\vert a\right\vert ^{i}\allowbreak r_{n-i}\left(
1|q\right) \allowbreak /\allowbreak \left( 1-q\right) ^{(n-i)/2}\allowbreak 
\\
&\leq &\allowbreak \left( 1-q\right) ^{-n/2}\sum_{i=0}^{n}\QATOPD[ ] {n}{i}%
q^{\binom{i}{2}}\left\vert a\sqrt{1-q}\right\vert ^{i}\allowbreak
r_{n-i}\left( 1|q\right) \allowbreak  \\
&\leq &\allowbreak \left( 1-q\right) ^{-n/2}\allowbreak r_{n}\left(
1|q\right) \allowbreak \left( -\left\vert a\right\vert \sqrt{1-q}\right)
_{n}\allowbreak  \\
&\leq &\allowbreak \left( -\left\vert a\right\vert \sqrt{1-q}\right)
_{\infty }\allowbreak \left( 1-q\right) ^{-n/2}\allowbreak r_{n}\left(
1|q\right) .
\end{eqnarray*}

v) An easy proof based on Carlitz formulae presented e.g. in Exercise
12.2(b) and 12.2(c) of \cite{IA} was given in \cite{AW}. It contained
however two misprints. Once concerning constant $\rho .$ Namely in all
estimates one should take $\left\vert \rho \right\vert $ naturally. The
second one concerned the following estimation: 
\begin{equation*}
w\left( x,y,\rho q^{k}|q\right) \leq (1-\rho ^{2}q^{2k})^{2}\allowbreak
+\allowbreak 2\sqrt{1-q}(1+\rho ^{2}q^{2k})|y\rho q^{k}|\allowbreak
+\allowbreak 4\rho ^{2}q^{2k}\allowbreak +\allowbreak (1-q)\rho
^{2}y^{2}q^{2k}
\end{equation*}
instead of 
\begin{equation*}
w\left( x,y,\rho q^{k}|q\right) \leq (1-\rho ^{2}q^{2k})^{2}\allowbreak
+\allowbreak 2(1-q)(1+\rho ^{2}q^{2k})|y\rho q^{k}|\allowbreak +\allowbreak
4\rho ^{2}q^{2k}\allowbreak +\allowbreak (1-q)\rho ^{2}y^{2}q^{2k}.
\end{equation*}
Now 
\begin{eqnarray*}
&&(1-\rho ^{2}q^{2k})^{2}\allowbreak +\allowbreak 2\sqrt{1-q}(1+\rho
^{2}q^{2k})|y\rho q^{k}|\allowbreak +\allowbreak 4\rho ^{2}q^{2k}\allowbreak
+\allowbreak (1-q)\rho ^{2}y^{2}q^{2k}\allowbreak  \\
&=&\allowbreak (1+\rho ^{2}q^{2k}+\sqrt{1-q}|y\rho q^{k}|)^{2}\leq
(1+\left\vert \rho \right\vert q^{k})^{4}
\end{eqnarray*}
since $\sqrt{1-q}|y|\leq 2$ for $y\in S\left( q\right) .$ Notice that $%
\prod_{k\geq 0}(1+\left\vert \rho \right\vert q^{k})^{4}\allowbreak
=\allowbreak (-\left\vert \rho \right\vert )_{\infty }^{4}.$ The last
inequality was suggested by the referee.
\end{proof}

\end{document}